\documentclass[12pt]{amsart}
\usepackage[colorlinks=true]{hyperref}
\usepackage{amsmath}
\usepackage[english,  activeacute]{babel}
\usepackage[utf8]{inputenc}
\usepackage{amssymb}
\usepackage{amsthm}
\usepackage{graphics,graphicx}
\usepackage{enumerate}
\usepackage{array}
\usepackage{bm}
\usepackage{a4wide}
\usepackage{color, url}
\usepackage{float}
\usepackage{multirow}
\newcommand\ChangeRT[1]{\noalign{\hrule height #1}}

\newcommand{\seqnum}[1]{\href{http://oeis.org/#1}{\underline{#1}}}

\theoremstyle{plain}
\newtheorem{theorem}{Theorem}[section]

\newtheorem{corollary}[theorem]{Corollary}

\theoremstyle{definition}

\newcommand{\C}{{\bm C}}
\newcommand{\ca}{{\bm c}}

\newcommand{\D}{{\bm D}}

\newcommand{\Cat}{{\mathcal{C}}}

\newcommand{\Ca}{{\bm C}}

\def\des{{\textsf{des}}}

\def\des{{\textsf{des}}}

\title{Descent distribution on Catalan words avoiding ordered pairs of Relations}
\date{\today}

\begin{document}

\author[J.-L. Baril]{Jean-Luc Baril}
\address{\noindent LIB, Université de Bourgogne, B.P. 47 870, 21078 DIJON-Cedex\\
France }
\email{barjl@u-bourgogne.fr}

\author[J. L. Ram\'{\i}rez]{Jos\'e L. Ram\'{\i}rez}
\address{Departamento de Matem\'aticas,  Universidad Nacional de Colombia,  Bogot\'a, Colombia}
\email{jlramirezr@unal.edu.co}

\begin{abstract}
 This work is a continuation of some recent articles presenting enumerative results for Catalan words avoiding one or a pair of  consecutive or classical patterns of length $3$. More precisely, we provide systematically the bivariate generating function for the number of Catalan words avoiding a given pair of relations with respect to the length and the number of descents. We also present several constructive bijections preserving the number of descents. As a  byproduct, we deduce the generating function for the total number of descents on all Catalan words of a given length and avoiding a pair of ordered relations. 
 \end{abstract}

\maketitle

\section{Introduction}

\begin{figure}[H]
\centering
\includegraphics[scale=0.8]{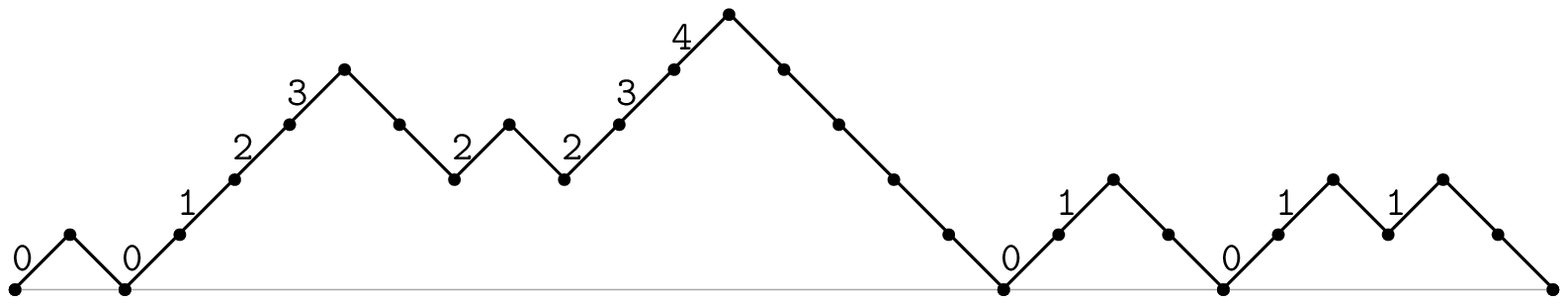}
\caption{Dyck path of the Catalan word \texttt{00123223401011}.} \label{fig1}
\end{figure}

A word $w=w_1w_2\cdots w_n$ over the set of non-negative integers is called a \emph{Catalan word}  if $w_1=\texttt{0}$ and $0\leq w_i\leq w_{i-1}+1$ for $i=2, \dots, n$. Let $\Cat_n$ denote the set of the Catalan words of length $n$. The cardinality of the set $\Cat_n$ is given by the Catalan number $C_n=\frac{1}{n+1}\binom{2n}{n}$, see  \cite[Exercise 80]{Stanley2}.   A \emph{Dyck path} of semilength $n$ is a lattice path of  $\mathbb{Z \times Z}$ running from $(0, 0)$ to $(2n, 0)$ that never passes below the $x$-axis and whose permitted steps are $U=(1, 1)$ and  $D=(1,-1)$.  For a Dyck path of semilength $n$, we associate a Catalan word in $\Cat_n$ formed  by the $y$-coordinate of each initial point of the up steps.  This construction is a bijection. For example, in Figure \ref{fig1}  we show the Dyck path associated to the Catalan word $\texttt{00123223401011}\in \Cat_{14}$.

 Catalan words have already been studied in the context of exhaustive generation
of Gray codes for growth-restricted words \cite{ManVaj}. More recently, Baril  et al. \cite{Baril, Baril2} study the distribution of descents on restricted Catalan words avoiding a pattern or a pair of patterns of length at most three. Ram\'irez and Rojas \cite{AlejaRam} also study the distribution of descents for Catalan words avoiding consecutive patterns of length at most three.  Baril, Gonz\'alez, and Ram\'irez \cite{BGR} enumerate Catalan words avoiding a classical pattern of length at most three according to the length and the value of the last symbol. They also give the exact value or an asymptotic for the expectation  of the last symbol. Also, we refer to  \cite{CallManRam, Toc}, where the authors study  several combinatorial statistics on the  polyominoes associated with words in $\Cat_n$. 
The goal of this work is to complement all these studies  by providing enumerative results for Catalan words avoiding a set of consecutive patterns defined from a pair of relations (see below for a formal definition),  with respect to the length and the number of descents.

The remaining of this paper is structured as follows. In Section 2, we introduce the notation that will be used in this work. In Section 3, we provide enumerative results for the number of Catalan words  avoiding a pattern defined by  an ordered pairs of relations. More precisely, in Section 3.1, we focus on the {\it constant cases}, i.e., pattern avoidances inducing a finite number of Catalan words independently of the length. All other sections handle exhaustively the remaining cases,  by providing bivariate generating functions with respect to the length and the number of descents. We also deduce the generating function with respect to the length for the total number of descents in Catalan words avoiding a given pattern.  
Below, Table 1 and 2 present an exhaustive list of the sequences counting Catalan words that avoid  an ordered pair of relations. Notice that the avoidance of an ordered pairs of relation can be equivalent to the avoidance of a consecutive pattern of length three (for instance, $(<,<)$  is equivalent to $\underline{012}$). In this case, the patterns were already studied in \cite{AlejaRam}.

\begin{table}[ht]
\centering
\begin{tabular}{|c|l|c|} \hline
$(X,Y)$ & Cardinality of $\Cat_n(X,Y)$, $n\geq 1$  & References \\ \hline \hline
$(\leq,\geq), (\leq,\neq)$ & 1, 2, 2, 2,  \dots  &  Section \ref{sec0}  \\ \hline
$(\leq,\leq)$ & 1, 2, 1, 2, 1, 2, \dots &  Section \ref{sec0}  \\ \hline
$(\neq,\leq)$ & 1, 2, 3, 3, 3,  \dots &  Section \ref{sec0}  \\ \hline
\end{tabular}
\caption{Number of Catalan words avoiding an ordered pair of relations: the constant cases.}
\label{Tab1}
\end{table}

\begin{table}[h]
\centering
\resizebox{0.96\textwidth}{!}{
\begin{tabular}{|c|c|c|c|} \hline
$(X,Y)$ & Cardinality of $\Cat_n(X,Y)$, $n\geq 1$ & OEIS & References\\ \hline \hline
$(=,=)$ & $\ca_{\underline{000}}(n)=\sum_{k=1}^{n}\binom{ k }{n-k}m_{k-1}$ & \seqnum{A247333} &  Theorem 2.8 of \cite{AlejaRam}, seq. $\ca_{\underline{000}}(n)$
 \\ \ChangeRT{1.1pt}

$(=,\geq)$ & 1, 2, 4, 10, 26, 72, 206, 606, 1820, 5558, \dots & \seqnum{A102407} &  Section \ref{sec1}  \\ \hline
$(\geq,=)$ & 1, 2, 4, 10, 26, 72, 206, 606, 1820, 5558, \dots & \seqnum{A102407} &  Section \ref{sec1} \\ 
\ChangeRT{1.1pt}

$(=,>)$ & $\C_{\underline{110}}(x)=\frac{1 - 2 x^2-\sqrt{1 - 4 x + 4 x^3}}{2 x (1 - x)}$  & \seqnum{A087626} & Theorem 2.6 of  \cite{AlejaRam}, seq. $\ca_{\underline{110}}(n)$    \\ \hline
$(>,=)$ & $\C_{\underline{100}}(x)=\frac{1 - 2 x^2-\sqrt{1 - 4 x + 4 x^3}}{2 x (1 - x)}$  & \seqnum{A087626}& Theorem 2.6 of  \cite{AlejaRam}, seq. $\ca_{\underline{100}}(n)$   \\ 
\ChangeRT{1.1pt}

$(=,\leq)$ & 1, 2, 3, 7, 17, 43, 114, 310, 861, 2433, \dots & \seqnum{A143013} &  Section \ref{sec2}    \\ \hline
$(\leq, =)$ & 1, 2, 3, 7, 17, 43, 114, 310, 861, 2433, \dots & \seqnum{A143013} & Section \ref{sec2}   \\ 
\ChangeRT{1.1pt}

$(=,<)$ & $\ca_{\underline{001}}(n)=\sum_{k=0}^{\lfloor (n-1)/2\rfloor} \frac{(-1)^k}{n - k}\binom{n - k}{k}\binom{2 n - 3 k}{n - 2 k-1}$ & \seqnum{A105633} & Theorem 2.3 of \cite{AlejaRam}, seq. $\ca_{\underline{001}}(n)$  \\ \hline
$(<,=)$ & $\ca_{\underline{011}}(n)=\sum_{k=0}^{\lfloor (n-1)/2\rfloor} \frac{(-1)^k}{n - k}\binom{n - k}{k}\binom{2 n - 3 k}{n - 2 k-1}$ & \seqnum{A105633} & Theorem 2.4 of \cite{AlejaRam}, seq. $\ca_{\underline{011}}(n)$  \\ \hline
$(<,>)$ & $\sum_{k=0}^{\lfloor (n-1)/2\rfloor} \frac{(-1)^k}{n - k}\binom{n - k}{k}\binom{2 n - 3 k}{n - 2 k-1}$ & \seqnum{A105633} & Section \ref{sec3}  \\ 
\ChangeRT{1.1pt}

$(=,\neq)$ &  1, 2, 4, 8, 17, 38, 89, 216, 539, 1374, \dots  & \seqnum{A086615} &  Section \ref{sec4}  \\ \hline
 $(\neq, =)$ &  1, 2, 4, 8, 17, 38, 89, 216, 539, 1374, \dots  & \seqnum{A086615} & Section \ref{sec4}   \\ 
\ChangeRT{1.1pt}

$(\geq,\geq)$ & $m_{n}$  (Motzkin numbers)   & \seqnum{A001006 } & Section \ref{sec5}   \\ \hline
$(<, <)$ &  $m_{n}$  (Motzkin numbers)  & \seqnum{A001006 } & Theorem 2.1 of \cite{AlejaRam}, seq. $\ca_n(\underline{012})$ \\ 
\ChangeRT{1.1pt}

$(\geq,>)$ & 1, 2, 5, 13, 35, 97, 275, 794, 2327, 6905, \dots& \seqnum{A082582}&  Section \ref{sec6}   \\ \hline
$(>,\geq)$ & 1, 2, 5, 13, 35, 97, 275, 794, 2327, 6905, \dots& \seqnum{A082582}& Section \ref{sec6}  \\ \hline
$(>,<)$ & 1, 2, 5, 13, 35, 97, 275, 794, 2327, 6905, \dots& \seqnum{A082582}& Section \ref{sec6}    \\ 
\ChangeRT{1.1pt}

$(\geq,\leq)$ &  $F_{n+1}$ (Fibonacci number) & \seqnum{A000045} &  Section \ref{sec7}  \\ \hline
$(\leq,<)$ &  $F_{n+1}$ (Fibonacci number) & \seqnum{A000045} & Section \ref{sec7}   \\ \hline
$(<,\leq)$ &  $F_{n+1}$ (Fibonacci number) & \seqnum{A000045} & Section \ref{sec7}   \\ 
\ChangeRT{1.1pt}

$(\geq,<)$ &  $2^{n-1}$ & \seqnum{A011782} &    Section \ref{sec8}  \\ \hline
$(\leq,>)$ &  $2^{n-1}$ & \seqnum{A011782} &    Section \ref{sec8}\\ 
\ChangeRT{1.1pt}

$(\geq,\neq)$ &  $\binom{n}{2}+1$ & \seqnum{A000124} &    Section \ref{sec8b} \\ 
\ChangeRT{1.1pt}

$(>,>)$ & $\ca_{\underline{210}}(n)=\sum_{k=0}^{\lfloor n/2\rfloor} \frac{1}{n - k}\binom{n - k}{k}\binom{n - k}{k+1}2^{n-2k-1}$ & \seqnum{A159771} & Theorem 2.9 of  \cite{AlejaRam}, seq. $\ca_{\underline{210}}(n)$  \\ 
\ChangeRT{1.1pt}

$(>,\leq)$ & $P_{n+1}$ (Pell numbers) & \seqnum{A000129} & Section \ref{sec9}.\\ 
\ChangeRT{1.1pt}

$(>,\neq)$ & 1, 2, 5, 13, 34, 90, 242, 660, 1821, 5073, \dots& New &     Section \ref{sec10} \\ 
\ChangeRT{1.1pt}

$(<,\geq)$ & $n$ & \seqnum{A000027} &    Section \ref{sec11}  \\ \hline

$(\neq,\geq)$ & $n$ & \seqnum{A000027} &     Section \ref{sec11}  \\ 
\ChangeRT{1.1pt}

$(<,\neq)$ &  1, 2, 3, 6, 12, 25, 54, 119, 267, 608, \dots& New &  Section \ref{sec12}   \\ 
\ChangeRT{1.1pt}

$(\neq,>)$ & 1, 2, 4, 9, 22, 56, 146, 388, 1048, 2869, \dots& \seqnum{A152225} &   Section \ref{sec13}   \\ 
\ChangeRT{1.1pt}

$(\neq,<)$ &  1, 2, 4, 8, 17, 37, 82, 185, 423, 978, \dots& \seqnum{A292460} &  Section \ref{sec14}  \\ 
\ChangeRT{1.1pt}

$(\neq,\neq)$ &  1, 2, 3, 6, 11, 22, 43, 87, 176, 362, \dots& \seqnum{A026418}& Section \ref{sec15}  \\ \hline
\end{tabular}}
\caption{Number of Catalan words avoiding an ordered pair of relations.}
\label{Tab2}
\end{table}%

\section{Notations}
 For an integer $r\geq 2$, a \emph{consecutive pattern} $p=\underline{p_1p_2\cdots p_r}$ is a word (of length $r$) over the set $\{ 0, 1, \dots, r-1\}$ satisfying the condition: if $j>0$ appears in $p$, then $j-1$ also appears in $p$. A Catalan word $w=w_1w_2\cdots w_n$  \emph{contains} the consecutive pattern $p=\underline{p_1p_2\cdots p_r}$ if  there exists a subsequence  $w_{i}w_{i+1}\cdots w_{i+r-1}$ (for some $i\geq 1$) of $w$  which is order-isomorphic to $p_1p_2\ldots p_r$. We say that $w$ \emph{avoids} the consecutive pattern $p$ whenever $w$ does not contain the consecutive pattern $p$. For example, the Catalan word $\texttt{0123455543}$ avoids the consecutive pattern $\underline{001}$ and contains one subsequence isomorphic to the pattern $\underline{210}$. More generally, we consider pattern $p$ as an ordered pair $p=(X,Y)$ of relations $X$ and $Y$ lying into the set $\{<,>,\leq,\geq,=,\neq\}$ (see for instance, Corteel et al. \cite{Cor}, Martinez and Savage  \cite{Sav}, and Auli and Elizalde \cite{AuliEli}). We will say that a Catalan word $w$ contains the pattern $p=(X,Y)$ if there exists $i\geq 1$ such that $w_i~X~w_{i+1}$ and $w_{i+1}~Y~w_{i+2}$. As an example, the pattern $(\neq, \geq)$ appears twice in the Catalan word $0123112$ on the triplets $231$ and $311$.  Notice that the avoidance of $(\neq, \geq)$ on Catalan words is equivalent to the avoidance of the four consecutive patterns $\underline{010}$, $\underline{011}$, $\underline{100}$, and $\underline{210}$. On the other hand, the avoidance of $(<,<)$ is equivalent to the consecutive pattern $\underline{012}$.

For $n\geq 0$ and for a given consecutive pattern $p$ or an ordered  pair of relations $p=(X,Y)$. Let  $\Cat_n(p)$ denote the set of Catalan words of length $n$ avoiding the consecutive pattern $p$. We denote by $\ca_{p}(n)$ the cardinality of  $\Cat_n(p)$, and we set $\Cat(p):=\bigcup_{n\geq 0}\Cat_n(p)$. We denote by $\des(w)$ the number of descents in $w$, i.e., the number of indices $i\geq 1$ such that $w_i>w_{i+1}$. Let $\Cat_{n,k}(p)$  denote the  set of Catalan words $w\in\Cat_n(p)$ such that $\des(w)=k$, and let $\ca_{p}(n,k):=|\Cat_{n,k}(p)|$.  Obviously, we have  $\ca_{p}(n)=\sum_{k=0}^{n-1}\ca_{p}(n,k)$. We introduce the bivariate generating function
$$\Ca_{p}(x,y):=\sum_{w\in \Cat(p)}x^{|w|}y^{\des(w)}=\sum_{n, k\geq 0}\ca_{p}(n,k)x^ny^k,$$
and we set 
$$\Ca_{p}(x):=\sum_{w\in\Cat(p)}x^{|w|}=\Ca_{p}(x,1).$$
The generating function for the total number of descents over all words in  $\Cat_n(p)$ is given by
$$\D_{p}(x):=\left.\frac{\partial \C_{p}(x,y)}{\partial y}\right|_{y=1}.$$

Throughout this work, we will often use the first return decomposition of a Catalan word $w$, which is $w=\texttt{0}(w'+1)w''$, where $w'$ and $w''$ are Catalan words, and where ($w'+1$) is the word obtained from $w'$ by adding 1 at all these symbols (for instance if $w'=\texttt{012012}$ then $(w'+1)=\texttt{123123}$). As an example, the first return decomposition of $w=\texttt{0122123011201}$ is given by setting $w'=\texttt{011012}$ and $w''=\texttt{011201}$.

\section{Enumeration}
Section 3.1 handles  pattern avoidances  inducing a finite number of Catalan words independently of the length. Section 3.1 to 3.17 handle exhaustively the other patterns.

\subsection{Constant cases}\label{sec0}
\leavevmode\par

The avoidance of $(\leq,\geq)$ on Catalan words is equivalent to the avoidance of $\underline{000}$, $\underline{010}$, $\underline{011}$, $\underline{110}$, and $\underline{120}$. So, Catalan words in $\Cat_n(\leq,\geq)$ are necessarily  of the form $\texttt{01}\cdots\texttt{n}$ or $\texttt{001}\cdots\texttt{(n-1)}$. Therefore $\ca_{(\leq,\geq)}(n)=2$ for all $n\geq 2$. 

The avoidance of $(\leq,\neq)$ on Catalan words is equivalent to the avoidance of $\underline{001}$, $\underline{010}$, $\underline{012}$, $\underline{110}$, and $\underline{120}$. So, we necessarily have   $\Cat_1(\leq,\neq)=\{\texttt{0}\}$ and $\Cat_n(\leq,\neq)=\{\texttt{0}^n,\texttt{0}\texttt{1}^{n-1} \}$  for $n\geq 2$. Therefore $\ca_{(\leq,\neq)}(n)=2$ for all $n\geq 2$.

The avoidance of $(\leq,\leq)$ on Catalan words is equivalent to the avoidance of $\underline{000}$, $\underline{001}$, $\underline{011}$, and $\underline{012}$. So, Catalan words in $\Cat_n(\leq,\leq)$ are words made up  of alternating zeros and ones except, possibly, for the last symbol, that is 
$$\Cat_n(\leq,\leq)=\begin{cases} (\texttt{01})^{n/2} \text{ or }(\texttt{01})^{(n-1)/2} \texttt{00},& \text{ if } n \text{ is even;}\\
(\texttt{01})^{(n-1)/2}\texttt{0},& \text{ if } n \text{ is odd.}\\
\end{cases}
$$
Therefore  we obtain $$\ca_{(\leq,\leq)}(n)=\begin{cases}
    2,& \text{ if } n \text{ is even};\\
      1,& \text{ if } n \text{ is odd}.
\end{cases}$$ 
Finally, the avoidance of $(\neq,\leq)$ on Catalan words is equivalent to the avoidance of $\underline{011}$, $\underline{012}$, $\underline{100}$, $\underline{101}$, and $\underline{201}$.  So, we have $\Cat_1(\neq,\leq)=\{\texttt{0}\}$, $\Cat_2(\neq,\leq)=\{\texttt{00}, \texttt{01}\}$, and for $n\geq 3$ 
$\Cat_n(\neq,\leq)=\{\texttt{0}^n, \texttt{0}^{n-1}\texttt{1}, \texttt{0}^{n-2}\texttt{10}\}$. Therefore, $$\ca_{(\neq,\leq)}(n)=\begin{cases}
    1,& \text{ if } n=1;\\
     2,& \text{ if } n=2;\\
     3,& \text{ if } n\geq 3.
\end{cases}$$ 
We refer to Table 1 for a summary of the results obtained in this part.

\subsection{Cases $\Cat(=, \geq)$ and $\Cat(\geq,=)$}\label{sec1}
\leavevmode\par

The avoidance of $(=,\geq)$ (resp. $(\geq,=)$) on Catalan words is equivalent to the avoidance of $\underline{000}$ and $\underline{110}$ (resp. $\underline{000}$ and $\underline{100}$). In \cite{BarilCol}, the authors established a bijection   between Catalan words avoiding \underline{100} and those avoiding \underline{110}: from left to right,  we replace each maximal factor $k^j(k-\ell)$, $j\geq 2$, $\ell\geq 1$, with the factor $k(k-\ell)^j$. This bijection preserves the number of descents and the avoidance of $\underline{000}$. 
For example, for $w = \texttt{0122123300}\in \Cat_{10}(\geq, =)$, we have the transformation
$$\texttt{01\textcolor{blue}{221}23300} \to \texttt{012112\textcolor{blue}{330}0} \to \texttt{0121123000}\in \Cat_{10}(=, \geq).$$
Therefore,  we necessarily have $\C_{(\geq, =)}(x,y)=\C_{(=, \geq)}(x,y)$, and below, we focus on the pattern $(=,\geq)$.
\begin{theorem}
We have
\begin{align*}
\C_{1}(x,y):=\C_{(=, \geq)}(x,y)&=\frac{1 - x - x^2 + 2 x y - \sqrt{(1 - x - x^2 + 2 x y)^2 - 4 x y (1 + x y)}}{2xy}.
\end{align*}
\end{theorem}
\begin{proof}
Let $w$ denote a non-empty Catalan word in $\Cat(=, \geq)=\Cat(\underline{000}, \underline{110})$, and let $w=\texttt{0}(w'+1)w''$ be the first return decomposition, where $w', w''\in \Cat(\underline{000}, \underline{110})$. If $w''=\epsilon$, then
$w=\texttt{0}(w'+1)$ with $w'$ possibly empty.  The generating function for this case is $x\C_{1}(x,y)$. If $w''$ is non-empty and $w'=\epsilon$, then  $w''$ cannot start with \texttt{00}. The corresponding generating function is $x A(x,y)$, where $A(x,y)$ is the bivariate generating function for non-empty Catalan words in $\Cat(\underline{000}, \underline{110})$ that do not start with  \texttt{00}. Counting using the complement we have
$$A(x,y)=\C_1(x,y) - 1 -
 \left[x^2\C_{1}(x,y) + x^2y(\C_1(x,y)-1-B(x,y))(\C_{1}(x,y)-1)\right],$$
where $B(x,y)$ is the bivariate generating function for Catalan words in $\Cat(=, \geq)$ that end with $aa$ for any $a\geq 0$ (we have to subtract this generating function to ensure the avoidance of $\underline{110}$). Since any word counted by $B(x,y)$ is obtained from a non-empty Catalan word in $\Cat(=, \geq)$ duplicating the last symbol unless the word ends with two repeated symbols, so this last  generating function satisfies $B(x,y)=x(\C_1(x,y)-1-B(x,y))$.

If $w'$ and $w''$ are non-empty, then  the generating function of this case is
$$E(x,y):=xy(\C_{1}(x,y)-1-B(x,y))(\C_1(x,y)-1).$$ Therefore, we have the functional equation
 \begin{align*}
\C_{1}(x,y) =1 + x\C_{1}(x,y)  + xA(x,y) +  E(x,y).
 \end{align*}
Solving this system of equations we obtain the desired result.
\end{proof}

 The  series expansion of the generating function $\C_{1}(x,y)$ is
\begin{align*}
1 + x + 2 x^2 + (3 + y) x^3 + (\bm{5} + \bm{5} y) x^4 + (8 + 16 y + 
    2 y^2) x^5 + (13 + 43 y + 16 y^2) x^6 + O(x^7).
\end{align*}
The Catalan words corresponding to the  bold coefficients in the above series are
\begin{align*}
\Cat_{4}(=, \geq)=\Cat_{4}(\underline{000}, \underline{110})=\{\texttt{00\textcolor{red}{10}}, \texttt{0011}, \texttt{0012}, \texttt{0\textcolor{red}{10}0},  \texttt{0\textcolor{red}{10}1}, \texttt{0112}, \texttt{01\textcolor{red}{20}},  \texttt{01\textcolor{red}{21}}, \texttt{0122}, \texttt{0123}\}.
\end{align*}

\begin{corollary} The g.f. for the cardinality of $\Cat(=, \geq)$  with respect to the length is
$$\C_{(=, \geq)}(x)=\frac{1 + x - x^2 - \sqrt{1 - 2 x - 5 x^2 - 2 x^3 + x^4}}{2 x}.$$
\end{corollary}
Using the  bijection given in Introduction between Catalan words and Dyck path, it is clear that the sequence  $\ca_{(=, \geq)}(n)=\ca_{\underline{000},\underline{100}}(n)$ also counts the 
number of all Dyck paths of semilength $n$ that avoid $DUDU$ (sequence \seqnum{A102407}). Then, using \cite{Sap} we have the combinatorial expression
$$\ca_{(=, \geq)}(n)=\ca_{\underline{000},\underline{110}}(n)=\sum_{j=0}^{\lfloor \frac{n}{2} \rfloor}\frac{1}{n  - j}\binom{n - j}{j}\sum_{i=0}^{n  - 2 j}\binom{n  - 2j}{i}\binom{j + i}{n - 2 j - i + 1},  \ n\geq 1.$$

\begin{corollary} The g.f. for the total number of descents on $\Cat(=, \geq)$  is
$$\D_{(=, \geq)}(x)=\frac{1 - 2 x - 3 x^2 + x^4 -(1-x-x^2) \sqrt{1 - 2 x - 5 x^2 - 2 x^3 + x^4}}{2 x\sqrt{1 - 2 x - 5 x^2 - 2 x^3 + x^4}}.$$
\end{corollary}
 The  series expansion of  $\D_{(=,\geq)}(x)$ is
\begin{align*}
 x^3 + 5 x^4 + 20 x^5 + 75 x^6 + 271 x^7+ 964 x^8+3397 x^9+O(x^{10})
\end{align*}
where the coefficient sequence does not appear in \cite{OEIS}.

\subsection{Cases $\Cat(=, \leq)$ and $\Cat(\leq,=)$}\label{sec2}
\leavevmode\par

The avoidance of $(=,\leq)$ (resp. $(\leq,=)$) on Catalan words is equivalent to the avoidance of $\underline{000}$ and $\underline{001}$ (resp. $\underline{000}$ and $\underline{011}$). In \cite{AlejaRam} (see Theorem 2.4), the authors established a bijection   between Catalan words avoiding \underline{011} and those avoiding \underline{001}.  From left to right,  they replace each factor $k^j(k+1)$ with the factor $k(k+1)^j$ $(j\geq 2)$. This bijection preserves the number of descents and the avoidance of $\underline{000}$. Therefore,  we necessarily have $\C_2(x,y):=\C_{(\leq, =)}(x,y)=\C_{(=, \leq)}(x,y)$.

\begin{theorem}
We have
\begin{align*}
\C_{2}(x,y):=\C_{(=, \leq)}(x,y)&=\frac{1 - x + 2 x y -  \sqrt{1 - 2 x + x^2 - 4 x^2 y - 4 x^3 y}}{2xy}.
\end{align*}
\end{theorem}
\begin{proof}
Let $w$ denote a non-empty Catalan word in $\Cat(=, \leq)=\Cat(\underline{000}, \underline{001})$, and let $w=\texttt{0}(w'+1)w''$ be the first return decomposition, where $w', w''\in \Cat(=, \leq)$. If $w''=\epsilon$, then
$w=\texttt{0}(w'+1)$ with $w'$ possibly empty.  The generating function for this case is $x\C_{2}(x,y)$. If $w''$ is non-empty and $w'=\epsilon$, then  $w''$ cannot start with \texttt{00} or \texttt{01}, so $w''=\texttt{0}$. The corresponding generating function is $x^2$. If $w'$ and $w''$ are non-empty, then  the generating function is
$xy(\C_{2}(x,y)-1)^2$. Therefore, we have the functional equation
 \begin{align*}
\C_{2}(x,y) =1 + x\C_{2}(x,y)  + x^2 +  xy(\C_{2}(x,y)-1)^2.
 \end{align*}
Solving this equation  we obtain the desired result.
\end{proof}

The  series expansion of the generating function $\C_{2}(x,y)$ is
\begin{align*}
1 + x + 2 x^2 + (2 + y) x^3 + (\bm{2} + \bm{5} y) x^4 + (2 + 13 y + 
    2 y^2) x^5 + (2 + 25 y + 16 y^2) x^6 + O(x^7).
\end{align*}
The Catalan words corresponding to the  bold coefficients in the above series are
\begin{align*}
\Cat_{4}(=, \leq)=\{\texttt{0\textcolor{red}{10}0}, \texttt{0\textcolor{red}{10}1}, \texttt{01\textcolor{red}{10}}, \texttt{01\textcolor{red}{20}}, \texttt{01\textcolor{red}{21}},  \texttt{0122}, \texttt{0123}\}.
\end{align*}
\begin{corollary} The g.f. for the cardinality of $\Cat(=, \leq)$ with respect to the length is
$$\C_{(=, \leq)}(x)=\frac{1 + x - \sqrt{1 - 2 x - 3 x^2 - 4 x^3}}{2 x}.$$
\end{corollary}
This generating function coincides with  those of the sequence \seqnum{A143013}, that is 
$$\ca_{(=, \leq)}(n)=\ca_{\underline{000},\underline{001}}(n)=\sum_{i=0}^{n}\sum_{k=1}^{n-i+1}\frac{1}{k}\binom{i - 1}{k - 1} \binom{k}{n - k - i + 1}\binom{k + i - 2}{i - 1},  \ n\geq 1.$$

Notice that sequence \seqnum{A143013} counts also the number of Motzkin paths with two kinds of level steps one of which is a final step.

\begin{corollary} The g.f. for the total number of descents on $\Cat(=, \leq)$  is
$$\D_{(=, \leq)}(x)=\frac{1 - 2 x - x^2 - 2 x^3 -(1-x) \sqrt{1 - 2 x - 3 x^2 - 4 x^3}}{2 x\sqrt{1 - 2 x - 3 x^2 - 4 x^3}}.$$
\end{corollary}

The  series expansion of  $\D_{(=,\leq)}(x)$ is
\begin{align*}
 x^3 + 5 x^4 + 17 x^5 + 57 x^6 + 188 x^7+ 610 x^8+1971 x^9+O(x^{10}),
\end{align*}
where the coefficient sequence does not appear in \cite{OEIS}.

\subsection{Cases $(=, <)$, $(<, =)$, and $(<,>)$}\label{sec3}
\leavevmode\par

Obviously, we have  $\Cat(=,<)=\Cat(\underline{001})$, $\Cat(<,=)=\Cat(\underline{011})$, and $\Cat(<,>)=\Cat(\underline{010},\underline{120})$. 
With the same bijection used at the beginning of Section \ref{sec2},   we have $\C_{(=,<)}(x,y)=\C_{(<,=)}(x,y)$, and we refer to \cite{AlejaRam} to see an expression of this generating function. On the other hand, a non-empty Catalan word in $\Cat(<,>)$ is either of the form ($i$) $\texttt{0}\alpha$ with $\alpha\in\Cat(<,>)$, ($ii$) $\texttt{0}(\alpha+1)$  with $\alpha\in\Cat(<,>)$, $\alpha\neq \epsilon$, or ($iii$) $\texttt{0}(\alpha+1)\beta$ where $\alpha$ ends with $a(a+1)$ and $\beta\in\Cat(<,>)$, $\beta\neq\epsilon$. We deduce the functional equation 
\begin{multline*}
\C_{(<,>)}(x)=1+x\C_{(<,>)}(x)+x(\C_{(<,>)}(x)-1)\\
+x(\C_{(<,>)}(x)-1)(\C_{(<,>)}(x)-1-x-x(\C_{(<,>)}(x)-1)),
\end{multline*}
which proves that $\C_{(<,>)}(x)=\C_{(=,<)}(x)=\C_{(<,=)}(x)$. Then, the sets $\Cat(\underline{011})=\Cat(<,=)$ and $\Ca(<,>)=\Ca(\underline{010},\underline{120})$ are in one-to-one correspondence, but the number of descents cannot be preserved (see for instance the list of Catalan words of length $3$). 

Below, we focus on the descent distribution over $\Cat(<,>)$.
\begin{theorem}
We have
\begin{align*}
\C_3(x,y):=\C_{(<,>)}(x,y)&=\frac{1 - 2 x + 2 x y - x^2 y -  \sqrt{1 - 4 x + 4 x^2 - 2 x^2 y + x^4 y^2}}{2 x y (1 - x)}.
\end{align*}
\end{theorem}

\begin{proof}
Let $w$ denote a non-empty Catalan word in $\Cat(<, >)=\Cat(\underline{010},  \underline{120})$.   Let $w=\texttt{0}(w'+1)w''$ be the first return decomposition, where $w', w''\in \Cat(<,>)$. If $w''=\epsilon$, then
$w=\texttt{0}(w'+1)$ with $w'$ possibly empty.  The generating function for this case is $x\C_{3}(x,y)$. If $w''$ is non-empty and $w'=\epsilon$, then  $w''$ is any word of $\Cat(<,>)$, and the corresponding generating function is $x (\C_3(x,y)-1)$. If $w'$ and $w''$ are non-empty, then  $w'\neq \texttt{0}$ or $w'$  does not end with an ascent $a(a+1)$, where $a\geq 0$. The generating function is
$$E(x,y):=xy(\C_{3}(x,y)-1-x-B(x,y))(\C_3(x,y)-1),$$ 
where $B(x,y)$ is the generating function for the Catalan words in $\Cat(<, >)$ ending with an ascent. Since any word counted by $B(x,y)$ is obtained from a non-empty Catalan word duplicating the last symbol plus 1, this   generating function is given by $B(x,y)=x(\C_3(x,y)-1)$. Therefore, we have the functional equation
 \begin{align*}
\C_{3}(x,y) =1 + x\C_{3}(x,y)+ x(\C_{3}(x,y)-1) + E(x,y).
 \end{align*}
Solving this system of equations we obtain the desired result.
\end{proof}

The  series expansion of the generating function $\C_{3}(x,y)$ is
\begin{align*}
1 + x + 2 x^2 + 
 4 x^3 + (\bm{8} + \bm{y}) x^4 + (16 + 6 y) x^5 + (32 + 24 y + y^2) x^6 + O(x^7).
\end{align*}
The Catalan words corresponding to the  bold coefficients in the above series are
\begin{align*}
\Cat_{4}(<,>)=\{\texttt{0000}, \texttt{0001}, \texttt{0011}, \texttt{0012}, \texttt{01\textcolor{red}{10}}, \texttt{0111}, \texttt{0112}, \texttt{0122},  \texttt{0123}\}.
\end{align*}
\begin{corollary} The g.f. for the cardinality of $\Cat(<, >)$ with respect to the length  is
$$\C_{(<,>)}(x)=\frac{1 -x^2 - \sqrt{1 - 4 x + 2 x^2 + x^4}}{2 (1-x)x}.$$
\end{corollary}
This generating function coincides with  the generating function of the sequence \seqnum{A105633}, and 
$$\ca_{(<,>)}(n)=\ca_{\underline{010},\underline{120}}(n)=\sum_{k=0}^{\lfloor (n-1)/2\rfloor} \frac{(-1)^k}{n - k}\binom{n - k}{k}\binom{2 n - 3 k}{n - 2 k-1},  \ n\geq 1.$$

Using the bijection given in Introduction between Catalan words and Dyck path, the sequence $\ca_{(<,=)}(n)=\ca_{\underline{011}}(n)$  counts the Dyck paths  of semilength $n$ avoiding $UUDU$ and $\ca_{(=,<)}(n)=\ca_{\underline{001}}(n)$ counts the Dyck paths of semilength $n$ avoiding $UDUU$ (cf. \cite{Sap}). Notice that the sequence \seqnum{A105633} also counts the number of Dyck paths of semilength $n+1$ with no pairs of consecutive valleys at the same height.  More generally,  it is proved in \cite{Sap2} that the generating function of the sequence
$$\ca_{\gamma=\underline{012\cdots k \ell(\ell+1)\cdots(\ell+s)}}(n), \quad 0\leq \ell \leq k \text{ and } s\geq 1,$$
satisfies the functional equation
$$\C_{\gamma}(x)=1+x\C_{\gamma}(x)^2-x^M\C_{\gamma}(x)^{M-(k+1-\ell)}\left(\C_{\gamma}(x) - \frac{1-(x\C_{\gamma}(x))^m}{1-x\C_{\gamma}(x)} \right),$$
where $M=\max\{k+1, s\}$ and $m=\min\{k+1,s\}$ (when $k=0$, $\ell=0$ and $s=1$ we retrieve the generating function $\C_{(=,<)}(x)=\C_{(<,>)}(x)$).

\begin{corollary} The g.f. for the total number of descents on $\Cat(<, >)$ is
$$\D_{(<,>)}(x)=\frac{1 - 4 x + 3 x^2  -(1-2x) \sqrt{1 - 4 x + 2 x^2 + x^4}}{2 (1 - x) x\sqrt{1 - 4 x + 2 x^2 + x^4}}.$$
\end{corollary}
The  series expansion of  $\D_{(<,>)}(x)$ is
\begin{align*}
 x^4 + 6 x^5 + 26 x^6 + 100 x^7+ 363 x^8+1277 x^9+O(x^{10}),
\end{align*}
where the coefficient sequence does not appear in \cite{OEIS}.

\subsection{Cases $\Cat(=, \neq)$ and $\Cat(\neq, =)$}\label{sec4}
\leavevmode\par
The avoidance of $(=,\neq)$ (resp.  $(\neq,=)$) on Catalan words is equivalent to the avoidance of $\underline{001}$ and $\underline{110}$ (resp.  $\underline{100}$ and $\underline{011}$). The following map is a bijection from $\Cat(=, \neq)$ to $\Cat(\neq, =)$ and it preserves the number of descents:   crossing $\sigma\in \Cat(=, \neq)$ from right to left,  we replace each factor $k(k+1)^j$ with the factor $k^j(k+1)$, $j\geq 2$, and we replace each factor $k(k-\ell)^j$ with the factor $k^j(k-\ell)$. Therefore,  we necessarily have $\C_{(=, \neq)}(x,y)=\C_{(\neq, =)}(x,y)$. 
 For example, for the word  $\texttt{01012323412300} \in \Cat_{14}(=, \neq)$, we have the transformation: 
\begin{multline*}
    \texttt{01012323412\textcolor{blue}{300}} \to \texttt{0101232341\textcolor{blue}{233}0} \to \texttt{010123234\textcolor{blue}{122}30} \\ \to \texttt{01012323\textcolor{blue}{411}230}      \to \texttt{0101232\textcolor{blue}{344}1230} \to \texttt{010123\textcolor{blue}{233}41230} \\ \to \texttt{01012\textcolor{blue}{322}341230} \to \texttt{0101\textcolor{blue}{233}2341230}
     \to \texttt{010\textcolor{blue}{122}32341230}  \to \texttt{01\textcolor{blue}{011}232341230}\\ \to \texttt{0\textcolor{blue}{100}1232341230} \to \texttt{\textcolor{blue}{011}01232341230}    \to \texttt{00101232341230} \in \Cat_{14}(\neq,=).
\end{multline*}

\begin{theorem}
We have
\begin{align*}
\C_{4}(x,y):=\C_{(=, \neq)}(x,y)&=\frac{1 - x + 2 x y - 2 x^2 y - \sqrt{1 - 2 x + x^2 - 4 x^2 y}}{2 x y (1 - x)}.
\end{align*}
\end{theorem}
\begin{proof}
Let $w$ denote a non-empty Catalan word in $\Cat(=, \neq)=\Cat(\underline{001}, \underline{110})$, and let $w=\texttt{0}(w'+1)w''$ be the first return decomposition, where $w', w''\in \Cat(=, \neq)$. If $w''=\epsilon$, then
$w=\texttt{0}(w'+1)$ with $w'$ possibly empty.  The generating function for this case is $x\C_{4}(x,y)$. If $w''$ is non-empty and $w'=\epsilon$, then  $w''=0^j$ for any $j\geq 1$ (to avoid \underline{001}). The corresponding generating function is $x (x/(1-x))$. If $w'$ and $w''$ are non-empty, then  the generating function is
$$E(x,y):=xy(\C_{4}(x,y)-1-B(x,y))(\C_4(x,y)-1),$$
where $B(x,y)$ is the bivariate generating function for Catalan words in $\Cat(=, \neq)$ that do not end with $aa$ ($a\geq 0$). It is clear that $B(x,y)=x(\C_4(x,y)-1)$. 
 Therefore, we have the functional equation
 \begin{align*}
\C_{4}(x,y) =1 + x\C_{4}(x,y)  + \frac{x^2}{1-x} +  E(x,y).
 \end{align*}
Solving this system of equations we obtain the desired result.
\end{proof}

 The  series expansion of the generating function $\C_{4}(x,y)$ is
\begin{align*}
1 + x + 2 x^2 + (3 + y) x^3 + (\bm{4} + \bm{4} y) x^4 + (5 + 10 y + 
    2 y^2) x^5 + (6 + 20 y + 12 y^2) x^6 + O(x^7).
\end{align*}
The Catalan words corresponding to the  bold coefficients in the above series are
\begin{align*}
\Cat_{4}(=, \neq)=\{\texttt{0000}, \texttt{0\textcolor{red}{10}0},  \texttt{0\textcolor{red}{10}1}, \texttt{0111},  \texttt{01\textcolor{red}{20}}, \texttt{01\textcolor{red}{21}}, \texttt{0122},  \texttt{0123}\}.
\end{align*}
\begin{corollary} The g.f. for the cardinality of $\Cat(=, \neq)$ with respect to the length is
$$\C_{(=, \neq)}(x)=\frac{1 + x - 2 x^2 - \sqrt{1 - 2 x - 3 x^2}}{2 (1 - x) x}.$$
\end{corollary}
This generating function coincides with  the generating function of the sequence of partial sums of the Motzkin numbers $m_n$ (see sequence \seqnum{A086615}), that is 
$$\ca_{(=, \neq)}(n)=\ca_{\underline{001},\underline{110}}(n)=\ca_{\underline{100},\underline{011}}(n)=\sum_{j=0}^{n-1} m_j,  \ n\geq 1,$$
where $m_n$ is the $n$-th Motzkin number. 

\begin{corollary} The g.f. for the total number of descents on $\Cat(=, \neq)$ is
$$\D_{(=, \neq)}(x)=\frac{1 - 2 x - x^2 -(1-x) \sqrt{1 - 2 x - 3 x^2}}{2 (1 - x) x\sqrt{1 - 2 x - 3 x^2}}.$$
\end{corollary}
The  series expansion of  $\D_{(=,\neq)}(x)$ is
\begin{align*}
 x^3 + 4 x^4 + 14 x^5 + 44 x^6+ 134 x^7+400 x^8+1184x^9+O(x^{10}),
\end{align*}
where the coefficient sequence corresponds to \seqnum{A097894} in \cite{OEIS}, which counts the  number of peaks at even height in all Motzkin paths of length $n+1$.

\subsection{Cases $\Cat(\geq, \geq)$ and $\Cat(<, <)$}\label{sec5}
\leavevmode\par

The avoidance of $(\geq,\geq)$ (resp.  $(<,<)$) on Catalan words is equivalent to the avoidance of $\underline{000}$, $\underline{100}$, $\underline{110}$, and $\underline{210}$ (resp.  $\underline{012}$). For the pattern $(<,<)=\underline{012}$, we refer to \cite{AlejaRam} to see an expression of $\C_{(<, <)}(x,y)$. Below, we prove that 
$$\C_{5}(x,y):=\C_{(\geq, \geq)}(x,y)=\C_{(<, <)}(x,y),$$ by exhibiting a bijection between $\Cat(\geq, \geq)$ and $\Cat(<, <)$ that preserves the length and the number of descents.
Let $w$ be a word in $\Cat(\geq, \geq)$,  we distinguish four cases:

$$\phi(w)=\left\{\begin{array}{ll}
\epsilon & \mbox{ if } w=\epsilon\\
\texttt{0} \phi(u) & \mbox{ if } w=\texttt{0}(u+1) \\
\texttt{01} (\phi(u)+1) & \mbox{ if } w=\texttt{00}(u+1)\\
\texttt{01}(\phi(v)+1)\phi(v') &  \mbox{ if } w=\texttt{0}(va+1)v'\\
\end{array}\right.,$$
where $u$, $va$,  and $v'$ are Catalan words in $\Cat(\geq, \geq)$, such that $u$ and $v$  are possibly empty, $v'$ is not empty, $v$ ends with $(a-1)$ when $v$ is not empty, and $v'$ does not start with $\texttt{00}$. Clearly, this map is a bijection from $\Cat(\geq, \geq)$ to $\Cat(<, <)$ that preserves the descent number. For instance, the image by $\phi$ of   $w=\texttt{0123010122} \in\Cat_{11}(\geq, \geq )$ is 
\begin{multline*}
\phi(w)=\texttt{01}\phi(\texttt{01})\cdot \phi(\texttt{010122})=
\texttt{01}\texttt{11}\cdot \texttt{01}\cdot \phi(\texttt{0122})\\=\texttt{011101}\texttt{0}\phi(\texttt{011})=\texttt{0111010001}\in\Cat(<, <).    
\end{multline*}
So, we deduce directly the following.

\begin{theorem}
We have
\begin{align*}
\C_{5}(x,y):=\C_{(\geq, \geq)}(x,y)&=\frac{1 - x - x^2 + x^2 y - \sqrt{(1 - x - x^2 (1 + y))^2 - 4 x^3 (1 + x) y}}{2 x^2 y}.
\end{align*}
\end{theorem}
 The  series expansion of the generating function $\C_{5}(x,y)$ is
\begin{multline*}
1 + x + 2 x^2 + (3 + y) x^3 + (\bm{5} + \bm{4} y) x^4 + (8 + 12 y + 
    y^2) x^5 + (13 + 31 y + 7 y^2) x^6 + O(x^7).
\end{multline*}
The Catalan words corresponding to the  bold coefficients in the above series are
\begin{align*}
\Cat_{5}(\geq, \geq)=\{ \texttt{00\textcolor{red}{10}}, \texttt{0011},  \texttt{0012},  \texttt{0\textcolor{red}{10}1}, \texttt{0112}, \texttt{01\textcolor{red}{20}}, \texttt{01\textcolor{red}{21}}, \texttt{0122}, \texttt{0123}\}.
\end{align*}
The coefficients of the  bivariate  generating function $\C_5(x,y)$ coincide with the array \seqnum{A114690}, which counts the number of Motzkin paths of length $n$ having $k$ weak ascents. 

\begin{corollary} The g.f. for the cardinality of $\Cat(\geq, \geq)$ with respect to the length is
$$\C_{(\geq, \geq)}(x)=\frac{1 - x - \sqrt{1 - 2 x - 3 x^2}}{2 x^2}.$$
\end{corollary}
This generating function coincides with  the generating function of the Motzkin numbers (sequence \seqnum{A001006})  $m_n=\sum_{k=0}^{\lfloor n/2 \rfloor}\binom{2n}{k}C_k$, where $C_k$ is the $k$-th Catalan number, that is 
$$\ca_{(\geq,\geq)}(n)=\ca_{\underline{000},\underline{100}, \underline{110}, \underline{210}}(n)=m_n,  \ n\geq 0.$$ 

\begin{corollary} The g.f. for the total number of descents on  $\Cat(\geq, \geq)$ is
$$\D_{(\geq, \geq)}(x)=\frac{1 - 2 x - 2 x^2 + x^3 -(1 - x - x^2) \sqrt{1 - 2 x - 3 x^2}}{2 x^2\sqrt{1 - 2 x - 3 x^2}}.$$
\end{corollary}
The series expansion of $\D_{(\geq, \geq)}(x)$ is 
$$x^3+ 4x^4+ 14x^5+45x^6+140x^7+427x^8+ 1288x^9+ O(x^{10}),$$
where the coefficients correspond to the sequence \seqnum{A005775} which counts triangular polyominoes of a given number of cells.

\subsection{Cases $\Cat(\geq, >)$, $\Cat(>,\geq)$, and $\Cat(>, <)$}\label{sec6}
\leavevmode\par
The avoidance of $\Cat(\geq, >)$ (resp. $\Cat(>,\geq)$, resp. $\Cat(>, <)$) on Catalan words is equivalent to the avoidance $\underline{110}$ and $\underline{210}$ (resp. $\underline{100}$ and $\underline{210}$, resp. $\underline{201}$ and $\underline{101}$).

From the bijection described in Section \ref{sec1}, the sets  $\Cat(\geq, >)=\Cat(\underline{110}, \underline{210})$ and \linebreak $\Cat(>,\geq)=\Cat(\underline{100},\underline{210})$ are in one-to-one correspondence (the bijection preserves also the descent number), which implies that  $\C_{(\geq, >)}(x,y)=\C_{(>,\geq)}(x,y)$.

On the other hand,  let us prove that 
$\Ca_{(\geq, >)}(x,y)=\Ca_{(>, <)}(x,y)$ by exhibiting a bijection $\psi$ preserving the descent number. Let $w$ be a word in $\Cat(\geq, >)$, then we distinguish five cases:
$$\psi(w)=\left\{\begin{array}{ll}
\epsilon & \mbox{ if } w=\epsilon\\
\texttt{0} \psi(u) & \mbox{ if } w=\texttt{0}u \\
\texttt{0} (\psi(u)+1) & \mbox{ if } w=\texttt{0}(u+1)\\
\texttt{0}(\psi(ua)+1)\texttt{0}\psi(v) &  \mbox{ if } w=\texttt{0}(u a(a+1)+1)v\\
\texttt{0}(\psi(v)+1)\texttt{0}& \mbox{ if } w=\texttt{01}v\\
\end{array}\right.,$$
where $u$, $ua(a+1)$,  and $v$ are Catalan words in $\Cat(\geq, >)$, such that $u$ is possibly empty and $v$ is not empty. Clearly, this map is a bijection from $\Cat(\geq, >)$ to $\Cat(>, <)$ that preserves the descent number.  For example, for the word  $\texttt{01234012343454} \in \Cat_{14}(\geq, >)$, we have the transformation: 
\begin{multline*}
    \texttt{0}\cdot  \texttt{1234} \cdot \texttt{012343454} \to \texttt{0}(\psi(\texttt{01}\cdot \texttt{2})+1)\texttt{0}\cdot \psi(\texttt{012343454}) \to \texttt{01}(\psi(\texttt{01})+2)\texttt{0}\cdot \texttt{0} (\psi(\texttt{01232343})+1)\\
\to   \texttt{01}\cdot\texttt{22}\cdot \texttt{00}\cdot \texttt{1}(\psi(\texttt{0121232})+2) \to    \texttt{0122001}\cdot \texttt{2}(\psi(\texttt{010121})+3) \\
\to  \texttt{01220012}\cdot \texttt{3}(\psi(\texttt{0121})+4)\texttt{3} \to  \texttt{012200123}\cdot \texttt{4}(\psi(\texttt{010})+5)\texttt{3} \\
\to \texttt{0122001234}\cdot \texttt{5}(\psi(\texttt{0})+6)5\cdot \texttt{3} \to \texttt{01220012345653} \in \Cat_{14}(>,<).
\end{multline*}

\begin{theorem}
We have
\begin{align*}
\C_{6}(x,y):=\C_{(\geq, >)}(x,y)&=\frac{1 - 2 x + x^2 y -  \sqrt{1 - 4 x + 4 x^2 - 2 x^2 y + x^4 y^2}}{2 x^2 y}.
\end{align*}
\end{theorem}
\begin{proof}
Let $w$ denote a non-empty Catalan word in $\Cat(\geq, >)=\Cat(\underline{110}, \underline{210})$, and let $w=\texttt{0}(w'+1)w''$ be the first return decomposition, where $w', w''\in \Cat(\geq, >)$. If $w''=\epsilon$, then
$w=\texttt{0}(w'+1)$ with $w'$ possibly empty.  The generating function for this case is $x\C_{6}(x,y)$. If $w''$ is non-empty and $w'=\epsilon$, then  $w''$ is any non-empty word in $\Cat(\geq, >)$, so the generating function is $x(\Ca_6(x,y)-1)$. If $w'$ and $w''$ are non-empty, then $w'=\texttt{0}$ or $w'$ has to finish with an ascent $a(a+1)$ ($a\geq 0$). Then the generating function is
$$E(x,y):=xy (x(\Ca_6(x,y)-1)) + xy (x(\Ca_6(x,y)-1))(\Ca_6(x,y)-1).$$ Therefore, we have the functional equation
 \begin{align*}
\C_{6}(x,y) =1 + x\C_{6}(x,y)  + x(\C_{6}(x,y)-1) +  E(x,y).
 \end{align*}
Solving this system of equations we obtain the desired result.
\end{proof}

 The  series expansion of the generating function $\C_{6}(x,y)$ is
\begin{align*}
1 + x + 2 x^2 + (4 + y) x^3 + (\bm{8} + \bm{5} y) x^4 + (16 + 18 y + 
    y^2) x^5 + (32 + 56 y + 9 y^2) x^6 + O(x^7).
\end{align*}
The Catalan words corresponding to the  bold coefficients in the above series are
\begin{align*}
\Cat_{4}(\geq, >)=\{ \texttt{0000}, \texttt{0001}, \texttt{00\textcolor{red}{10}}, \texttt{0011}, \texttt{0012}, 
\texttt{0\textcolor{red}{10}0}, \texttt{0\textcolor{red}{10}1},  \texttt{0111},  \texttt{0112}, \texttt{01\textcolor{red}{20}}, \texttt{01\textcolor{red}{21}}, \texttt{0122}, \texttt{0123}\}.
\end{align*}
The coefficient sequence of the  bivariate generating function $\C_6(x,y)$ coincides with \seqnum{A273717}, which counts  the number of bargraphs of semiperimeter $n$ having $k$ $L$-shaped corners $(n\geq 2, k\geq 0)$.

\begin{corollary} The g.f. for the cardinality of $\Cat(\geq, >)$ with respect to the length is
$$\C_{(\geq, >)}(x)=\frac{1 - 2 x + x^2 - \sqrt{1 - 4 x + 2 x^2 + x^4}}{2 x^2}.$$
\end{corollary}
This generating function coincides with  the generating function of the sequence \seqnum{A082582}, that is 
$$\ca_{(\geq, >)}(n)=\ca_{\underline{110},\underline{210}}(n)=\sum_{k=0}^{n}\sum_{j=0}^{n-k}  \frac{1}{j+1}\binom{n - k - 1}{j}\binom{k}{j}\binom{k + j + 2}{j},  \ n\geq 1.$$
Notice that this sequence counts also Dyck paths of semilength $n+1$ avoiding $UUDD$.

\begin{corollary} The g.f. for the total number of descents on $\Cat(\geq, >)$ is
$$\D_{(\geq, >)}(x)=\frac{1 - 4 x + 3 x^2 - (1 - 2 x) \sqrt{1 - 4 x + 2 x^2 + x^4}}{2 x^2\sqrt{1 - 4 x + 2 x^2 + x^4}}.$$
\end{corollary}
The series expansion of $\D_{(\geq, >)}(x)$ is 
$$x^3+5x^4+20x^5+74x^6+263x^7+914x^8+3134x^9+O(x^{10}),$$
where the coefficients correspond to the sequence \seqnum{A273718}, which counts also the total number of descents in all bargraphs of semiperimeter $n-1$.

\subsection{Cases $\Cat(\geq, \leq)$, $\Cat(\leq,<)$, and $\Cat(<, \leq)$}\label{sec7}
\leavevmode\par

The avoidance of $(\geq, \leq)$ (resp. $(\leq,<)$, resp. $(<, \leq)$) on Catalan words is equivalent to the avoidance of 
$\underline{000}$, $\underline{001}$, $\underline{100}$, $\underline{101}$, and  $\underline{201}$ (resp. $\underline{001}$ and $\underline{012}$, resp. $\underline{011}$ and $\underline{012}$).

From the bijection described in Section \ref{sec2},   $\Cat(\leq, <)=\Cat(\underline{001}, \underline{012})$ and $\Cat(<,\leq)=\Cat(\underline{011}, \underline{012})$ are in one-to-one correspondence by preserving the descent number, which implies that $\C_{(\leq,<)}(x,y)=\C_{(<,\leq)}(x,y)$.

\begin{theorem}
We have
\begin{align*}
\C_{7}(x,y):=\C_{(\leq,<)}(x,y)&=\frac{1 + x^2 - x^2y}{1 - x - x^2y}.
\end{align*}
\end{theorem}
\begin{proof}
Let $w$ denote a non-empty Catalan word in $\Cat(\leq, <)=\Cat(\underline{001}, \underline{012})$. From the conditions we have the decomposition  $\texttt{0}\texttt{0}^j$, $\texttt{01}\texttt{1}^j$ or $\texttt{01}\texttt{1}^jw'$, where $j\geq 0$ and  $w'$ is a non-empty word in $\Cat(\leq, <)$. Therefore, we have the functional equation
 \begin{align*}
\C_{7}(x,y) =1 + \frac{x}{1-x}  + \frac{x^2}{1-x} + \frac{x^2y}{1-x}\C_7(x,y).
 \end{align*}
Solving this equation we obtain the desired result.
\end{proof}

 The  series expansion of the generating function $\C_{7}(x,y)$ is
\begin{align*}
1 + x + 2 x^2 + (2 + y) x^3 + (\bm{2} + \bm{3} y) x^4 + (2 + 5 y + 
    y^2) x^5 + (2 + 7 y + 4 y^2) x^6 + O(x^7).
\end{align*}
The Catalan words corresponding to the  bold coefficients in the above series are
\begin{align*}
\Cat_{4}(\leq, <)=\{ \texttt{0000}, \texttt{0\textcolor{red}{10}0}, \texttt{0\textcolor{red}{10}1}, \texttt{01\textcolor{red}{10}}, \texttt{0111}\}.
\end{align*}
The coefficients of the generating function $\C_{7}(x,y)$  coincide with the table \seqnum{A129710}.

\begin{corollary} The g.f. for the cardinality of  $\Cat(\leq, <)$ with respect to the length is
$$\C_{(\leq, <)}(x)=\frac{1}{1-x-x^2}.$$
\end{corollary}
This generating function coincides with that of the Fibonacci sequence $F_{n+1}$, that is 
$$\ca_{(\leq, <)}(n)=\ca_{\underline{001},\underline{012}}(n)=F_{n+1},  \ n\geq 0.$$

 \begin{corollary} The g.f. for the total number of descents on  $\Cat(\leq, <)$ is
$$\D_{(\leq, <)}(x)=\frac{x^3 (1 + x)}{(1 - x - x^2)^2}.$$
\end{corollary}
The series expansion is 
$$x^3+3x^4+7x^5+ 15x^6 +30x^7 +58x^8 +109x^9+O(x^{10}),$$
where the coefficient sequence corresponds to \seqnum{A023610}. 
\medskip

Finally, we can prove easily that 
$\Ca_{(\geq, \leq)}(x,y)=\Ca_{(\leq, <)}(x,y)$ by exhibiting a bijection $\phi$ preserving the descent number. Let $w$ be a word in $\Cat(\leq, <)$, then we distinguish three cases: ($i$) $\phi(\epsilon)=\epsilon$; ($ii$) if $w=\texttt{0}^j$, $j\geq 1$, then we set $\phi(w)=\texttt{01}^{j-1}$; and ($iii$) if $w=\texttt{0}\texttt{1}^j$, $j\geq 1$, then $\phi(w)=\texttt{012}\cdots (j-1)(j-1)$;  ($iv$) if $w=\texttt{01}^jw'$, $j\geq 1$, then $\phi(w)=\texttt{01}\cdots (j-1)(\phi(w')+j)\texttt{0}$.
For example, for the word  $\texttt{01101101111011} \in \Cat_{14}(\leq, <)$, we have the transformation: 
\begin{multline*}
\phi(\texttt{01}^2\texttt{01101111011})=\texttt{01}(\phi(\texttt{01}^2\texttt{01111011})+2)\texttt{0}\\
=\texttt{0123}(\phi(\texttt{01}^4\texttt{011})+4)\texttt{20}=\texttt{01234567}(\phi(\texttt{01}^2)+8)\texttt{420}\\=\texttt{01234557899}(\phi(\epsilon)+8)\texttt{420}=\texttt{01234567899420}\in\Cat_{14}(\geq,\leq).\end{multline*}

\subsection{Cases $\Cat(\geq, <)$ and $\Cat(\leq, >)$}\label{sec8}
\leavevmode\par
The avoidance of $\Cat(\geq, <)$ (resp. $\Cat(\leq, >)$) on Catalan words is equivalent to the avoidance of $\underline{001}$, $\underline{101}$, and $\underline{201}$ (resp. $\underline{010}$, $\underline{110}$. and $\underline{120}$). Below, we will see  that for all $n\geq 1$, the cardinalities of $\Cat_n(\geq, <)$ resp. $\Cat_n(\leq, >)$ are the same, but the distributions of the number of descents do not coincide.

\begin{theorem}
We have
\begin{align*}
\C_{8}(x,y):=\C_{(\geq, <)}(x,y)&=\frac {1-x+x^2-x^2y}{1-2x+x^2-x^2y}.
\end{align*}
\end{theorem}
\begin{proof}
Let $w$ denote a non-empty Catalan word in $\Cat(\geq, <)=\Cat(\underline{101}, \underline{110}, \underline{201})$ and let $w=\texttt{0}(w'+1)w''$ be the first return decomposition, where $w', w''\in \Cat(\geq, <)$. If $w'$ is empty, then we have $w=0^j$, $j\geq 0$. Otherwise,  we have $w=\texttt{0}(w'+1)0^j$, $j\geq 0$.  The generating function $\C_8(x,y)$ satisfies the functional equation 
$$\C_8(x,y)=\frac{1}{1-x}+ x(\C_{8}(x,y)-1)+\frac{x^2y}{1-x}(\C_{8}(x,y)-1),$$
which gives the result.
\end{proof}

 The  series expansion of the generating function $\C_{8}(x,y)$ is
\begin{align*}
1 + x + 2 x^2 + (3 + y) x^3 + (\bm{4} + \bm{4} y) x^4 + (5 + 10 y + 
    y^2) x^5 + (6 + 20 y + 6 y^2) x^6 + O(x^7).
\end{align*}
The Catalan words corresponding to the  bold coefficients in the above series are
\begin{align*}
\Cat_{4}(\geq, <)=\{ \texttt{0000}, \texttt{01\textcolor{red}{10}},  \texttt{0111},\texttt{0\textcolor{red}{10}0}, \texttt{01\textcolor{red}{20}}, \texttt{01\textcolor{red}{21}}, \texttt{0122}, \texttt{0123}\}.
\end{align*}

The coefficient sequence of the  bivariate generating function $\C_8(x,y)$ coincides with \seqnum{A034867}, then 
$$[x^ny^k]\C_8(x,y)=\binom{n}{2k+1}.$$

\begin{corollary} The g.f. for the cardinality of  $\Cat(\geq, <)$ with respect to the length is
$$\C_{(\geq, <)}(x)=\frac{1-x}{1-2x},$$
where the $n$th term is $2^{n-1}$.
\end{corollary}

 \begin{corollary} The g.f. for the total number of descents  on  $\Cat(\geq, <)$  is
$$\D_{(\geq, <)}(x)=\frac{x^3}{(1 -2x)^2}.$$
\end{corollary}
The series expansion is 
$$x^3+4x^4+12x^5+ 32x^6 +80x^7 +192x^8 +448x^9+O(x^{10}),$$
where the coefficient sequence corresponds to \seqnum{A001787}, i.e., the $n$-th term is $(n-2)2^{n-3}$. 

\begin{theorem}
We have
\begin{align*}
\C_{(\leq, >)}(x,y)&=\frac {1-x}{1-2x}.
\end{align*}
\end{theorem}
\begin{proof}
Let $w$ denote a non-empty Catalan word in $\Cat(\leq, >)=\Cat(\underline{010}, \underline{110}, \underline{120})$ and let $w=\texttt{0}(w'+1)w''$ be the first return decomposition, where $w', w''\in \Cat(\geq, <)$. If $w'$ is empty, then we have $w=0w''$. Otherwise,  we have $w=\texttt{0}(w'+1)$.  The generating function $\C_{(\leq, >)}(x,y)$ satisfies the functional equation 
$$\C_{(\leq, >)}(x,y)=1+ x\C_{(\leq, >)}(x,y) + x(\C_{(\leq, >)}(x,y)-1),$$
which gives the result.
\end{proof}

\subsection{Case $\Cat(\geq , \neq)$}\label{sec8b}
\leavevmode\par
The avoidance of $\Cat(\geq, \neq)$ on Catalan words is equivalent to the avoidance of $\underline{001}$, $\underline{102}$, $\underline{201}$, $\underline{110}$, and $\underline{210}$. This means that any Catalan word of length $n$ in $\Cat(\geq, \neq)$  is either of the form $w=01\ldots (k-1)k^{n-k}$, $k\geq 0$, or  $w=01\ldots (k-1)m^{n-k}$ where $0\leq m<k-1$ and $k\geq 1$. Therefore, we can deduce easily that 
$$ \C_9(x,y):=\C_{(\geq, \neq)}(x,y)=1+\frac{x}{(1-x)^2}+\frac{yx^3}{(1-x)^3},$$
which implies the following.

\begin{theorem}
We have
\begin{align*}
\C_{(\geq, \neq)}(x,y)&=\frac {1-2x+2x^2-x^3+x^3y}{\left( 1-x \right) ^{3}}.
\end{align*}
\end{theorem}

 The  series expansion of the generating function $\C_{(\geq, \neq)}(x,y)$ is
\begin{align*}
1 + x + 2 x^2 + (3 + y) x^3 + (\bm{4} + \bm{3} y) x^4 + (5 + 6 y) x^5 + (6 + 10 y) x^6 + O(x^7).
\end{align*}
The Catalan words corresponding to the  bold coefficients in the above series are
\begin{align*}
\Cat_{4}(\geq, \neq)=\{ \texttt{0000}, \texttt{0\textcolor{red}{10}0},  \texttt{0111},\texttt{01\textcolor{red}{20}}, \texttt{01\textcolor{red}{21}},  \texttt{0122}, \texttt{0123}\}.
\end{align*}

\begin{corollary} The g.f. for the cardinality of  $\Cat(\geq , \neq)$ with respect to the length is
$$\C_{(\geq, \neq)}(x)=\frac {1-2x+2x^2}{ \left( 1-x \right) ^{3}},$$
where the $n$th term is $1+\binom{n}{2}$.
\end{corollary}

 \begin{corollary} The g.f. for the total number of descents  on  $\Cat(\geq, \neq)$ is
$$\D_{(\geq, \neq)}(x)=\frac{x^3}{(1 -x)^3}.$$
\end{corollary}
The series expansion is 
$$x^3+3x^4+6x^5+ 10x^6 +15x^7 +21x^8 +28x^9+O(x^{10}),$$
where the coefficient sequence corresponds to \seqnum{A000217}, i.e., the $n$-th term is $\binom{n-1}{2}$.

\subsection{Case $\Cat(>, \leq)$}\label{sec9}
\leavevmode\par

The avoidance of $\Cat(>, \leq)$ on Catalan words is equivalent to the avoidance of $\underline{100}$, $\underline{101}$, and $\underline{201}$.

\begin{theorem}
We have
\begin{align*}
\C_{10}(x,y):=\C_{(>, \leq)}(x,y)&=\frac{1 - x - x^2 y}{1 - 2 x - x^2 y}.
\end{align*}
\end{theorem}
\begin{proof}
Let $w$ denote a non-empty Catalan word in $\Cat(>, \leq)=\Cat(\underline{100}, \underline{101}, \underline{201})$, and let $w=\texttt{0}(w'+1)w''$ be the first return decomposition, where $w', w''\in \Cat(>, \leq)$. If $w''=\epsilon$, then
$w=\texttt{0}(w'+1)$ with $w'$ possibly empty.  The generating function for this case is $x\C_{10}(x,y)$. If $w''$ is non-empty and $w'=\epsilon$, then  $w''$ is any non-empty word in $\Cat(>, \leq)$, so the generating function is $x(\Ca_{10}(x,y)-1)$. If $w'$ and $w''$ are non-empty, then $w'$ does not start with the prefix $\texttt{01}$ or $\texttt{00}$, otherwise $w$ would contain the pattern $\underline{100}$ or $\underline{201}$, respectively. Therefore, $w''=\texttt{0}$ and the generating function is $x^2y(\Ca_{10}(x,y)-1)$.
Summarizing, we have the functional equation
 \begin{align*}
\C_{10}(x,y) =1 + x\C_{10}(x,y)  + x(\C_{10}(x,y)-1) +  x^2y(\Ca_{10}(x,y)-1).
 \end{align*}
Solving this equation we obtain the desired result.
\end{proof}

The  series expansion of the generating function $\C_{10}(x,y)$ is
\begin{align*}
1 + x + 2 x^2 + (4 + y) x^3 + (\bm{8} + \bm{4} y) x^4 + (16 + 12 y + 
    y^2) x^5 + (32 + 32 y + 6 y^2) x^6 + O(x^7).
\end{align*}
The Catalan words corresponding to the  bold coefficients in the above series are
\begin{align*}
\Cat_{4}(>, \leq)=\{\texttt{0000}, \texttt{0001}, \texttt{00\textcolor{red}{10}}, \texttt{0011}, \texttt{0012}, \texttt{01\textcolor{red}{10}}, \texttt{0111}, \texttt{0112}, \texttt{01\textcolor{red}{20}}, \texttt{01\textcolor{red}{21}}, \texttt{0122}, \texttt{0123}\}.
\end{align*}
The coefficients of the above series coincide with the array \seqnum{A207538}. Notice that this array coincides with the coefficients of the Pell polynomials.

\begin{corollary} The g.f. for the cardinality of $\Cat(\leq, <)$ with respect to the length  is
$$\C_{(>, \leq)}(x)=\frac{1 - x - x^2}{1 - 2 x - x^2}.$$
\end{corollary}

This generating function coincides with  the generating function of the Pell numbers $P_{n+1}$. The Pell sequence is defined by $P_n=2P_{n-1} + P_{n-2}$ for $n\geq 2$, with the initial values $P_0=0$ and $P_1=1$ (see sequence \seqnum{A000129}).

\begin{corollary} The g.f. for the total number of descents on $\Cat(>, \leq)$  is
$$\D_{(>, \leq)}(x)=\frac{x^3}{(1 - 2x - x^2)^2}.$$
\end{corollary}
The series expansion of $\D_{(>, \leq)}(x)$ is 
$$x^3+4x^4+14x^5+44x^6+131x^7+376x^8+1052x^9+O(x^{10}),$$
where the coefficients correspond to the sequence \seqnum{A006645}.

\subsection{Case $\Cat(>, \neq)$}\label{sec10}
\leavevmode\par
The avoidance of $\Cat(>, \neq)$ on Catalan words is equivalent to the avoidance of $\underline{101}$, $\underline{201}$, and  $\underline{210}$, which means that $\Cat(>, \neq)$ is the set of Catalan words without double successive descents and without valleys.
Therefore, such a nonempty word is of the form $\texttt{0}u$,  or $\texttt{0}(u+1)$, or 
$\texttt{0}(u+1)u'\texttt{0}v$ where $u,v\in\Cat(>, \neq)$ and $u'=k(k+1)\cdots \ell$ with $\ell\geq k\geq 1$ and if $u+1$ ends with $a$ then we set $k=a$, and if $u$ is empty then we set $k=1$. We deduce the following  functional equation for $\C_{11}(x,y):=\C_{(>, \neq)}(x,y)$.
$$\C_{11}(x,y)=1+x\C_{11}(x,y)+x(\C_{11}(x,y)-1)+yx^2\frac{x}{1-x}\C_{11}(x,y)^2.$$

Then, we have the following.

\begin{theorem}
We have
\begin{align*}
\C_{11}(x,y):=\C_{(>, \neq)}(x,y)&=\frac{\left(1-2x - \sqrt {1 - 4 x + 4 x^2 - 4 x^3 y}
 \right)  \left(1-x \right)}{2x^3y}.
\end{align*}
\end{theorem}
The series expansion of the generating function $\C_{11}(x,y)$ is
$$1+x+2x^2+(4+y)x^3+(\bm{8}+\bm{5}y)x^4+(16+18y)x^5+(32+56y+2y^2)x^6+ O(x^7).$$
The Catalan words corresponding to the bold coefficients in the above series are 
$$\Cat_4(>, \neq)=\{\texttt{0000}, \texttt{0001}, \texttt{00\textcolor{red}{10}}, \texttt{0011}, \texttt{0012}, \texttt{0\textcolor{red}{10}0}, \texttt{01\textcolor{red}{10}}, \texttt{0111}, \texttt{0112},  \texttt{01\textcolor{red}{20}},\texttt{01\textcolor{red}{21}}, \texttt{0122}, \texttt{0123}\}$$

\begin{corollary} The g.f. for the cardinality of $\Cat(>, \neq)$ with respect to the length  is
$$\C_{(>, \neq)}(x)=\frac{\left(1-2x - \sqrt {1 - 4 x + 4 x^2 - 4 x^3}
 \right)  \left(1-x \right)}{2x^3}.$$
\end{corollary}

The coefficient sequence of the series expansion does not appear in \cite{OEIS}.

\begin{corollary} The g.f. for the total number of descents  on $\Cat(>, \neq)$ is
$$\D_{(>, \neq)}(x)=
\frac {\left(1-x \right)\left(1-4x+4x^2-2x^3-(1-2x)\sqrt {1 - 4 x + 4 x^2 - 4 x^3} \right) }{2x^3\sqrt {1 - 4 x + 4 x^2 - 4 x^3}}.$$
\end{corollary}
The series expansion of $\D_{(>, \neq)}(x)$ is 
$$x^3+5x^4+18x^5+60x^6+196x^7+632x^8+2015x^9+O(x^{10}),$$
where the coefficient sequence does not appear in \cite{OEIS}. 

\subsection{Cases $\Cat(<, \geq)$ and $\Cat(\neq, \geq)$ }\label{sec11}
\leavevmode\par
The avoidance of $\Cat(<, \geq)$  (resp. $\Cat(\neq, \geq)$) on Catalan words is equivalent to the avoidance of $\underline{010}$, $\underline{120}$, and  $\underline{011}$ (resp. $\underline{100}$, $\underline{011}$, $\underline{210}$, $\underline{010}$, and $\underline{120}$). 
A nonempty Catalan word in $\Cat(<, \geq)$  is  of the form $\texttt{0}^k(u+1)$, $k\geq 1$ where $u$ is either empty or  $u=\texttt{012} \ldots k$, $k\geq 1$. Moreover, it is easy to check that $\Cat(<, \geq)=\Cat(\neq, \geq)$.
Then, we deduce the following.

\begin{theorem}
We have
\begin{align*}
\C_{12}(x,y):=\C_{(<, \geq)}(x,y)=\C_{(\neq, \geq)}(x,y)=1+\frac{1}{1-x}\frac{x}{1-x}={\frac {1-x+x^2}{ \left(1-x \right) ^{2}}}.
\end{align*}
\end{theorem}
The $n$th term of the series expansion is $n$.

\subsection{Case $\Cat(<, \neq)$}\label{sec12}
\leavevmode\par
The avoidance of $\Cat(<, \neq)$ on Catalan words is equivalent to the avoidance of $\underline{010}$, $\underline{012}$, and  $\underline{120}$. We set $\C_{13}(x,y):= \C_{(<, \neq)}(x,y)$.

\begin{theorem}
We have
\begin{align*}
\C_{13}(x,y)=\frac {1 - 2 x^2 - x^3 + 2 x^2 y-(1+x)\sqrt {1 - 2 x - x^2 + 2 x^3 + x^4 - 4 x^3 y}}{2{x}^{2}y}.
\end{align*}
\end{theorem}

\begin{proof}
Let $w$ denote a non-empty Catalan word in $\Cat(<, \neq)=\Cat(\underline{010}, \underline{012}, \underline{120})$, and let $w=\texttt{0}(w'+1)w''$ be the first return decomposition, where $w', w''\in \Cat(<, \neq)$. If $w'$ is empty, then the generating function for these words is $x\C_{13}(x,y)$. If $w''$ is empty, then $w=\texttt{01}(u+1)$, where $u\in  \Cat(<, \neq)$, and the generating function for these words is $x^2\C_{13}(x,y)$. If $w=\texttt{011}u$, where $u$ is a nonempty word in $\Cat(<, \neq)$, then the generating function for these words is $x^3y(\C_{13}(x,y)-1)$. If $w=\texttt{01}(u+1)v$, where $u,v\in \Cat(<, \neq)$, $u$ ending with $aa$, $a\geq 0$, and $v$ nonempty, then the generating function for these words is $x^3y(\C_{13}(x,y)-1)^2$. If  $w=\texttt{01}(u+1)v$, where $u,v\in \Cat(<, \neq)$, $u$ ending with a descent and $v$ nonempty, then  the generating function for these words is $x^2yB(x,y)(\C_{13}(x,y)-1)$, where $B(x,y)$ is the generating function for Catalan words in $\Cat(<, \neq)$ ending by a descent. By considering the complement, we obtain easily that $$B(x,y)=\C_{13}(x,y)-1-x-xB(x,y)-x(\C_{13}(x,y)-1)-x^2\C_{13}(x,y).$$
Summarizing, we have the following functional equation 
\begin{multline*}
\C_{13}(x,y)=1+x\C_{13}(x,y)+x^2\C_{13}(x,y)+x^3y(\C_{13}(x,y)-1)\\+x^3y(\C_{13}(x,y)-1)^2+x^2yB(x,y)(\C_{13}(x,y)-1),
\end{multline*}
which gives the result.
\end{proof}
The series expansion of the generating function $\C_{13}(x,y)$ is
$$1 + x + 2 x^2 + 3 x^3 + (5 + y) x^4 + (\bm{8 + 4 y}) x^5 + (13 + 12 y) x^6 +O(x^7)$$
The Catalan words corresponding to the bold coefficients in the above series are 
\begin{multline*}
    \Cat_5(<, \neq)=\left\{\texttt{00000}, \texttt{00001}, \texttt{00011}, \texttt{001\textcolor{red}{10}}, \texttt{00111},\texttt{00112},\right.\\
    \left. \texttt{01\textcolor{red}{10}0}, \texttt{01\textcolor{red}{10}1}, \texttt{011\textcolor{red}{10}},  \texttt{01111}, \texttt{01112},\texttt{01122}\right\}\end{multline*}
\begin{corollary} The g.f. for the cardinality of $\Cat(<, \neq)$ with respect to the length  is
$$\C_{(<, \neq)}(x)={\frac {1-{x}^{3}-(1+x)\sqrt {1 - 2 x - x^2 - 2 x^3 + x^4}}{2{x}^{2}}}
.$$
\end{corollary}
The coefficient sequence of the series expansion does not appear in \cite{OEIS}.
\begin{corollary} The g.f. for the total number of descents on $\Cat(<, \neq)$ is
$$\D_{(<, \neq)}(x)=
\frac {1 - 4 x^2 - 4 x^3 + 2 x^5 - x^6-(1+x)(1-2x^2-x^3)\sqrt { 1 - 2 x - x^2 - 2 x^3 + x^4}}{2x^2(1 + x) \sqrt { 1 - 2 x - x^2 - 2 x^3 + x^4}}.$$
\end{corollary}
The series expansion of $\D_{(<, \neq)}(x)$ is 
$$x^4+4x^5+12x^6+35x^7+97x^8+262x^9+O(x^{10}),$$
where the coefficient sequence does not appear in \cite{OEIS}.

\subsection{Case $\Cat(\neq , >)$}\label{sec13}
\leavevmode\par
The avoidance of $\Cat(\neq,>)$ on Catalan words is equivalent to the avoidance of $\underline{010}$, $\underline{210}$, and  $\underline{120}$. 
\begin{theorem}
We have
\begin{align*}
\C_{14}(x,y):=\C_{(\neq , >)}(x,y)&=\frac{1 - 2 x + 2 x^2 y - \sqrt{1 - 4 x + 4 x^2 - 4 x^3 y}}{2 x^2 y}.
\end{align*}
\end{theorem}
\begin{proof}
Let $w$ denote a non-empty Catalan word in $\Cat(\neq , >)=\Cat(\underline{010}, \underline{120}, \underline{210})$, and let $w=\texttt{0}(w'+1)w''$ be the first return decomposition, where $w', w''\in \Cat(\neq , >)$. If $w''=\epsilon$, then
$w=\texttt{0}(w'+1)$ with $w'$ possibly empty.  The generating function for this case is $x\C_{14}(x,y)$. If $w''$ is non-empty and $w'=\epsilon$, then  $w''$ is any non-empty word in $\Cat(\neq,>)$, so the generating function is $x(\Ca_{14}(x,y)-1)$. If $w'$ and $w''$ are non-empty, then $w'\neq \texttt{0}$ (to avoid $\underline{010}$),  $w'$ does not end with an  ascent $a(a+1)$ ($a\geq 0$) (to avoid $\underline{120}$) or  $w'$ does not end with a  descent $ab$ ($a>b$) (to avoid $\underline{210}$). Then, $w'$ ends with $aa$ ($a\geq 0$) and the generating function is $x^2y(\Ca_{14}(x,y)-1)^2$.
Summarizing, we have the functional equation
 \begin{align*}
\C_{14}(x,y) =1 + x\C_{14}(x,y)  + x(\C_{14}(x,y)-1) +  x^2y(\Ca_{14}(x,y)-1)^2.
 \end{align*}
Solving this equation  we obtain the desired result.
\end{proof}

The  series expansion of the generating function $\C_{14}(x,y)$ is
\begin{align*}
1 + x + 2 x^2 + 4 x^3 + (\bm{8 + y}) x^4 + (16 + 6 y) x^5 + (32 + 24 y) x^6 + O(x^7).
\end{align*}
The Catalan words corresponding to the  bold coefficients in the above series are
\begin{align*}
\Cat_{4}(\neq , >)=\{\texttt{0000}, \texttt{0001}, \texttt{0011}, \texttt{0012}, \texttt{01\textcolor{red}{10}}, \texttt{0111}, \texttt{0112}, \texttt{0122}, \texttt{0123}\}.
\end{align*}
\begin{corollary} The g.f. for the cardinality of  $\Cat(\neq, >)$ with respect to the length  is
$$\C_{(\neq, >)}(x)=\frac{1 - 2 x + 2 x^2 - \sqrt{1 - 4 x + 4 x^2 - 4 x^3}}{2x^2}.$$
\end{corollary}
This generating function coincides with  the generating function of the sequence \seqnum{A152225} that counts Dyck paths of a given length with no peaks at height $0$ (mod 3) and no valleys at height $2$ (mod 3). 

\begin{corollary} The g.f. for the total number of descents on $\Cat(\neq, >)$  is
$$\D_{(\neq, >)}(x)=\frac{1 - 4 x + 4 x^2 - 2 x^3 - (1 - 2 x)\sqrt{1 - 4 x + 4 x^2 - 4 x^3}}{2x^2\sqrt{1 - 4 x + 4 x^2 - 4 x^3}}.$$
\end{corollary}
The series expansion of $\D_{(\neq,>)}(x)$ is 
$$x^4+6x^5+24x^6+84x^7+280x^8+912x^9+O(x^{10}),$$
where the coefficient sequence does not appear in \cite{OEIS}. 

\subsection{Case $\Cat(\neq , <)$}\label{sec14}
\leavevmode\par

The avoidance of $\Cat(\neq, <)$ on Catalan words is equivalent to the avoidance of $\underline{012}$, $\underline{101}$, $\underline{102}$, and $\underline{201}$. 

\begin{theorem}
We have
\begin{align*}
\C_{15}(x,y):=\C_{(\neq , <)}(x,y)&=\frac{1 - x - x^2 - \sqrt{(1 - x - x^2)^2 - 4 x^3 y}}{2 x^3 y}.
\end{align*}
\end{theorem}
\begin{proof}
Let $w$ denote a non-empty Catalan word in $\Cat(\neq , <)=\Cat(\underline{012}, \underline{101}, \underline{201})$, and let $w=\texttt{0}(w'+1)w''$ be the first return decomposition, where $w', w''\in \Cat(\neq , <)$. If $w''=\epsilon$, then
$w=\texttt{0}(w'+1)$ with $w'$ possibly empty.  For $w'\neq \epsilon$ and to avoid \underline{012},  we have  $w'=0\texttt{1}(w'''+1)$, where $w'''\in \Cat(\neq , <)$. The generating function for this case is $x+x^2\C_{15}(x,y)$. If $w''$ is non-empty and $w'=\epsilon$, then  $w''$ is any non-empty word in $\Cat(\neq,<)$, so the generating function is $x(\Ca_{15}(x,y)-1)$. If $w'$ and $w''$ are non-empty, then $w'=\texttt{0}u$ and $w''=\texttt{0}v$, where $u, v\in \Cat(\neq , <)$. Then the generating function for this case is  $x^3y\C_{15}(x,y)$.
Summarizing, we have the functional equation
 \begin{align*}
\C_{15}(x,y) =1 + x+x^2\C_{15}(x,y)  + x(\Ca_{15}(x,y)-1) + x^3y\C_{15}(x,y).
 \end{align*}
Solving this equation we obtain the desired result.
\end{proof}

The  series expansion of the generating function $\C_{15}(x,y)$ is
\begin{align*}
1 + x + 2 x^2 + (3 + y) x^3 + (\bm{5 + 3 y}) x^4 + (8 + 9 y) x^5 + (13 + 
    22 y + 2 y^2) x^6 + O(x^7).
\end{align*}
The Catalan words corresponding to the  bold coefficients in the above series are
\begin{align*}
\Cat_{4}(\neq , >)=\{\texttt{0000}, \texttt{0001}, \texttt{00\textcolor{red}{10}}, \texttt{0011}, \texttt{0\textcolor{red}{10}0}, \texttt{01\textcolor{red}{10}}, \texttt{0111}, \texttt{0112}\}.
\end{align*}
The coefficients of the generating function $\C_{15}(x,y)$  coincide with the array \seqnum{A114711}. 
\begin{corollary} The g.f. for the cardinality of  $\Cat(\neq, >)$ with respect to the length is
$$\C_{(\neq, <)}(x)=\frac{1 - x - x^2 - \sqrt{1 - 2 x - x^2 - 2 x^3 + x^4}}{2x^3}.$$
\end{corollary}
This generating function coincides with  the generating function of the sequence \seqnum{A292460} that gives the number of $U_{k}D$-equivalence classes of \L{}ukasiewicz paths (see \cite{BKP}), which is a shift of the sequences   \seqnum{A292460}, \seqnum{A004148}, and   \seqnum{A203019}. Then we have 
$$\ca_{(\neq, <)}(n)=\ca_{\underline{012},\underline{101},\underline{201}}(n)=\sum_{k=0}^{\lfloor \frac{n+1}{2}\rfloor}(-1)^k\binom{n - k - 1}{k}m_{n + 1 - 2 k},  \ n\geq 0.$$

\begin{corollary} The g.f. for the total number of descents on   $\Cat(\neq, <)$  is
$$\D_{(\neq, <)}(x)=\frac{1 - 2 x - x^2 + x^4 - (1 - x - x^2)\sqrt{1 - 2 x - x^2 - 2 x^3 + x^4}}{2x^3\sqrt{1 - 2 x - x^2 - 2 x^3 + x^4}}.$$
\end{corollary}
The series expansion of $\D_{(\neq,<)}(x)$ is 
$$x^3+3x^4+9x^5+26x^6+71x^7+191x^8+508x^9+O(x^{10}),$$
where the coefficient sequence does not appear in \cite{OEIS}.

\subsection{Case $\Cat(\neq, \neq)$}\label{sec15}
\leavevmode\par
The avoidance of $\Cat(\neq, \neq)$ on Catalan words is equivalent to the avoidance of $\underline{010}$, $\underline{012}$, $\underline{101}$, $\underline{201}$, and  $\underline{120}$. 
\begin{theorem}
We have
\begin{align*}
\C_{16}(x,y):=\C_{(\neq, \neq)}(x,y)&=\frac{1 - x - x^2 - \sqrt{(1 - x - x^2)^2 - 4 x^4 y}}{2 x^4 y}.
\end{align*}
\end{theorem}
\begin{proof}
Let $w$ denote a non-empty Catalan word in $\Cat(\neq, \neq)=\Cat(\underline{010}, \underline{012}, \underline{101}, \underline{120}, \underline{201}, \underline{210})$, and let $w=\texttt{0}(w'+1)w''$ be the first return decomposition, where $w', w''\in \Cat(\neq, \neq)$. If $w''=\epsilon$, then
$w=\texttt{0}(w'+1)$ with $w'$ possibly empty.  For $w'\neq \epsilon$ and to avoid \underline{012},  we have  $w'=0\texttt{1}(w'''+1)$, where $w'''\in \Cat(\neq , \neq)$. The generating function for this case is $x+x^2\C_{16}(x,y)$. If $w''$ is non-empty and $w'=\epsilon$, then  the generating function is $x(\Ca_{16}(x,y)-1)$. If $w'$ and $w''$ are non-empty, then $w'$ is a non-empty word in $\Cat(\neq, \neq)$ such that the last two symbols are equal (to avoid \underline{210} and \underline{120}) and $w''=\texttt{0}w'''$, where $w'''\in\Cat(\neq, \neq)$ (to avoid \underline{101} and \underline{201}). Then the generating function for this case is  
$$E(x,y):=xy (x^2\Ca_{16}(x,y))(x \Ca_{16}(x,y))).$$
Summarizing, we have the functional equation
 \begin{align*}
\C_{16}(x,y) =1 + x+x^2\C_{16}(x,y)  + x(\Ca_{16}(x,y)-1) + E(x,y).
 \end{align*}
Solving this equation we obtain the desired result.
\end{proof}

Notice that $\C_{16}(x,y)=\C_{15}(x,xy)$. The  series expansion of the generating function $\C_{16}(x,y)$ is
\begin{align*}
1 + x + 2 x^2 + 3 x^3 + (\bm{5 + y}) x^4 + (8 + 3 y) x^5 + (13 + 9 y) x^6 + O(x^7).
\end{align*}
The Catalan words corresponding to the  bold coefficients in the above series are
\begin{align*}
\Cat_{4}(\neq , \neq)=\{\texttt{0000}, \texttt{0001}, \texttt{0011}, \texttt{01\textcolor{red}{10}}, \texttt{0111}, \texttt{0112}\}.
\end{align*}
\begin{corollary} The g.f. for the cardinality of  $\Cat(\neq, \neq)$ with respect to the length is
$$\C_{(\neq, \neq)}(x)=\frac{1 - x - x^2 - \sqrt{1 - 2 x - x^2 + 2 x^3 - 3 x^4}}{2x^4}.$$
\end{corollary}
This generating function coincides with  the generating function of the sequence \seqnum{A026418} that counts  ordered trees with a given number of  edges and having no branches of length~$1$. Then, we have 
$$\ca_{(\neq, \neq)}(n)=\ca_{\underline{010}, \underline{012}, \underline{101}, \underline{120}, \underline{201}, \underline{210}}(n)=\sum_{k=0}^{\lfloor \frac{n+1}{2}\rfloor}\binom{n - k}{k}m_{k},  \ n\geq 1.$$

\begin{corollary} The g.f. for the total number of descents on  $\Cat(\neq, \neq)$ is
$$\D_{(\neq, \neq)}(x)=\frac{1 - 2 x - x^2 + 2 x^3 - x^4 - (1 - x - x^2)\sqrt{1 - 2 x - x^2 + 2 x^3 - 3 x^4}}{2x^4\sqrt{1 - 2 x - x^2 + 2 x^3 - 3 x^4}}.$$
\end{corollary}
The series expansion of $\D_{(\neq,<)}(x)$ is 
$$x^4+3x^5+9x^6+22x^7+55x^8+131x^9+O(x^{10}),$$
where the coefficient sequence does not appear in \cite{OEIS}.

\end{document}